\documentclass[12pt]{article}
\usepackage{amsfonts}
\usepackage{amsmath}
\usepackage{epsfig}

\title{Continued Fractions and the $4$-Color Theorem}
\author{Richard Evan Schwartz \thanks{Supported by N.S.F. Grant DMS-2102802}}

\newtheorem{theorem}{Theorem}[section]

\newtheorem{lemma}[theorem]{Lemma}

\newtheorem{conjecture}[theorem]{Conjecture}

\def\startproof{{\bf {\medskip}{\noindent}Proof: }}

\def\endproof{$\spadesuit$  \newline}

\def\C{\mbox{\boldmath{$C$}}}%
\def\Q{\mbox{\boldmath{$Q$}}}%
\def\R{\mbox{\boldmath{$R$}}}%
\def\Z{\mbox{\boldmath{$Z$}}}%

\begin{document}
\maketitle

\begin{abstract}
  We study the geometry of some proper
  $4$-colorings of the vertices of sphere triangulations
  with degree sequence
  $6,...,6,2,2,2$.
  Such triangulations are the simplest examples
  which have non-negative combinatorial curvature.
  The examples we construct, which are roughly
  extremal in some sense, are based on a novel
  geometric interpretation of continued fractions.
  We will also present a conjectural sharp
  ``isoperimetric inequality'' for colorings
of this kind of triangulation.
\end{abstract}

\section{Introduction}

\subsection{Background}

The Four Color Theorem, first proved (with the
assistance of a computer) by
Wolfgang Haken and Kenneth Appel in $1976$, is one of the most
famous results in mathematics.  See
[{\bf W\/}] for a thorough discussion.
Here is
one formulation.  If you have any triangulation
of the $2$-sphere, it is possible to color the
vertices using $4$ colors such that no two
adjacent vertices have the same coloring.
This is called a {\it proper vertex $4$-coloring\/}.

Certainly one can properly $4$-color the vertices of
a tetrahedron.  A proper vertex $4$-coloring of a
triangulation $\cal Z$ (with the same colors)
canonically defines a simplicial
map $f$ from the sphere to the tetrahedron:  Just
map each vertex of $\cal Z$ to the
like-colored vertex of the tetrahedron and then
extend linearly to the faces.  The map $f$ in turn
defines a $2$-coloring of the faces of $\cal Z$.
One colors a face of $\cal Z$ black if $f$ is
orientation preserving on that face, and otherwise
white.

The associated face $2$-coloring has the property that
around each vertex
the number of black faces is congruent mod $3$ to
the number of white faces.  This derives from the
property that $3$ triangles of the tetrahedron meet
around each vertex. We call a face
 $2$-coloring with this property a {\it good coloring\/}.
 Conversely, a good coloring for
$\cal Z$ defines a simplicial map to the
tetrahedron and thus a proper $4$-coloring of the vertices of $\cal Z$.
So, an equivalent formulation of the
$4$-color theorem is that every triangulation
of the sphere has a good coloring.

So far, the Four Color Theorem only has computer-assisted
proofs.  Perhaps one can get insight into the result by
looking at examples of good colorings.
The good coloring version has a geometric feel to it, and so perhaps
some geometric insight might help.
The purpose of this paper is to look at the
geometry of these good colorings in some
special cases.

A {\it triangulation of non-negative combinatorial curvature\/} is
one in which the maximum degree is $6$.
All the vertices have degree $6$ except for a list
$v_1,...,v_k$ which have degrees $2 \leq d_1,...,d_k<6$.
Euler's Formula gives the condition on the degrees:
\begin{equation}
  \label{GB}
  \sum_{i=1}^k (6-d_i)=12.
  \end{equation}
In particular $k \leq 12$.
The quantity
$$\frac{\pi}{3} \times (6-d_i)$$
is the {\it combinatorial curvature\/} at $v_i$.
Equation \ref{GB} translates into a discrete version of the
Gauss-Bonnet theorem, which says that the
total combinatorial curvature is $4\pi$.

The triangulations of non-negative combinatorial
curvature form an attractive family to study.
In [{\bf T\/}], William Thurston organized these triangulations
into moduli spaces.
To give some idea of how this works, a triangulation of
non-negative combinatorial curvature defines a flat
cone structure on the sphere with non-negative
curvature: we just make all the
triangles unit equilateral triangles. The set of all
triangulations with the same list $d_1,..,d_k$ of defects includes
in the moduli space of flat cone structures on spheres
with appropriately prescribed singularities.  So, even though the
triangulations don't exactly vary continuously, one
can think of them as special points inside moduli
spaces consisting of structures which do vary
continuously.  Also, if the triangulations involve many triangles
and the defects are well spread out, one can imagine
that the defects almost vary continuously.

Given the nice structure of the totality
of such triangulations, it seems like an interesting idea
to study the space of good  colorings as a kind of
partially defined bundle over these moduli spaces.  Perhaps
the structure of such colorings is related somehow to
the placement of the defects.  As the defects vary around,
perhaps the good colorings vary in a nice way to some extent.
I imagine that the total picture, seen all at once,
would be spectacularly beautiful.
All this is very speculative. In spite of making a lot of
computer experiments over the years -- every
time I teach the graph theory class at Brown I play with this
project -- I don't have much to report.
In this very modest paper I will consider
the simplest cases.  The cases I have in mind are where
$k=3$ and $d_1=d_2=d_3=2$. 

\subsection{The Continued Fraction Colorings}

The $6,...,6,2,2,2$ trianglulations are
indexed by the nonzero Eisenstein integers. An {\it Eisenstein
  integer\/} is a number of the form
$$a+b \alpha, \hskip 30 pt
a,b \in \Z, \hskip 15 pt \alpha=\frac{1+i \sqrt 3}{2}.$$
(It is more common to use $\omega=\alpha^2$ in place of
$\alpha$, but $\alpha$ is more convenient for us.)
To see the connection, let $\cal E$ denote the ring of
Eisenstein integers. The points of $\cal E$ are naturally
the vertices of an equilateral triangulation $\cal T$ of $\C$.
Given some nonzero $\beta \in \cal E$, we let
$\beta \cal T$ be the bigger equilateral triangulation
obtained by multiplying the whole picture by $\beta$.
We let $G_{\beta}$ denote the
group of symmetries generated by order $3$ rotations
in the vertices of $\beta \cal T$. The quotient
$${\cal T\/}(\beta) = {\cal T\/}/G_{\beta}$$
is the desired triangulation.

Since all the vertices of ${\cal T\/}(\beta)$ have even
degree, ${\cal T\/}(\beta)$ always has a good coloring.
One just colors the triangles alternately black and white
in a checkerboard pattern.  Indeed, this face coloring
corresponds to a proper $3$-coloring of the vertices.

Figure 1 shows a different good coloring for
$\beta=2+3\alpha$.  We call this coloring ${\cal C\/}(2+3\alpha)$.
To get the triangulation of the sphere, cut out the specially outlined central
rhombus, fold it up like a taco, and glue the edges together in
pairs.   Figure 1 is
really showing part of the orbifold universal cover of the coloring.
Figures 7 and 8 below show more elaborate examples.

\begin{center}
  \resizebox{!}{3in}{\includegraphics{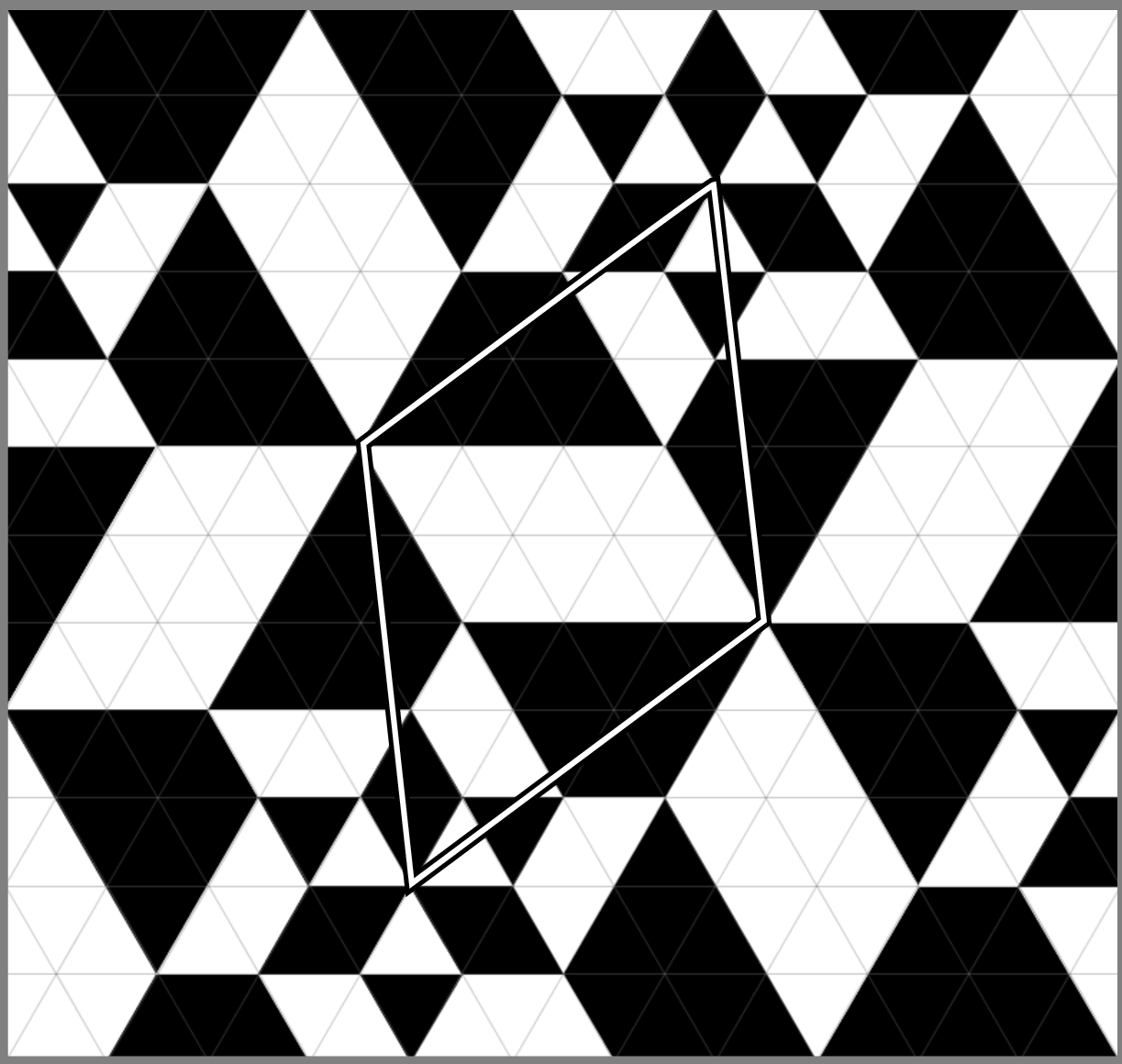}}
\newline
{\bf Figure 1:\/} The good coloring
${\cal C\/}(2+3\alpha)$ of
${\cal T\/}(2+3\alpha)$.
  \end{center}

  The coloring in Figure 1 is (at least experimentally)
  extremal in a certain sense.  Define the
  {\it fold count\/} of a good coloring to be the
  number of edges which form black-white interfaces.
  For the alternating coloring, the fold count is
  $3/2$ times the number of faces.  For ${\cal T\/}(2 + 3 \alpha)$
  this comes to $57$.   In contrast, ${\cal C\/}(2+3\alpha)$
  has fold count $23$.  It seems that ${\cal C\/}(2+3\alpha)$
 minimizes the fold count amongst all
 good colorings of ${\cal T\/}(2+3\alpha)$.

  The number of good colorings of ${\cal T\/}(\beta)$ grows
  exponentially with $|\beta|$, for an easy reason.
  In any good coloring, if we can find a vertex of degree
  $6$ where the colors alternate, we can switch the colors
  and get another good coloring.  In terms of the original
  vertex coloring, the neighbors of such a vertex $v$ are
  colored using just $2$ colors, and so we have an option
  to switch the color of $v$ to the other available color.
  Starting with the alternating
  coloring, we can take a large family of non-adjacent
  vertices and make these swaps according to any
  binary sequence we like.
  On the other hand, these examples seem rather similar to
  the alternating coloring.
  Their fold count is linear in the number of the faces.
  
  The colorings with smaller fold counts, and in particular
  with minimal fold counts, seem to be much more rigid and
  interesting.  To use an analogy from statistical
  mechanics, the colorings with large fold count are
  sort of like a fluid or a gas, and the colorings with
  small fold count are more like solid crystals.
  The example in Figure 1 is part of an infinite sequence of
  examples.  We call these examples ${\cal C\/}(\beta)$, where
  $\beta$ ranges over the primitive Eisenstein integers.
  (An Eisenstein integer $\beta=a+b \alpha$ {\it primitive\/} if it is not an
  integer multiple of another one. Equivalently, $a$ and $b$ are
  relatively prime.)
  We will see that ${\cal C\/}(\beta)$ is a special good coloring
  of  ${\cal T\/}(\beta)$ which geometrically implements
  the continued fraction expansion of $a/b$.   We call
  these colorings {\it continued fraction colorings\/}.
  
  \subsection{Properties of the Continued Fraction Colorings}
  
  The main result of this paper is that these continued fraction
  colorings exist, but I will prove some additional results about them.
  One interesting property is that these colorings have the same
  number of black and white triangles. Indeed, any good coloring
  of a triangulation having all even degrees has the same number
  of black and white triangles.  See \S \ref{even} for
  the quick proof which I learned from Kasra Rafi.

Our next result
  concerns the asymptotics of the fold
  count for the continued fraction colorings.
  Let $\{a_n/b_n\} \in (0,1)$ be a sequence of rationals.
    Let $${\cal T\/}_n={\cal T\/}(a_n+b_n \alpha), \hskip 30 pt
    {\cal C\/}_n={\cal C\/}(a_n+b_n \alpha).$$
    Let $f_n$ denote the fold count of ${\cal C\/}_n$ and
    let $F_n$ denote the number of faces in ${\cal T\/}_n$.

    \begin{theorem}
      \label{ASY}
      The following is true about the continued fraction colorings.
      \begin{enumerate}
        
    \item If $\{a_n/b_n\}$ converges to an irrational limit then
      $\lim_{n \to \infty} f_n/F_n = 0$.

  \item If $\{a_n/b_n\}$ is the sequence of continued fraction
    approximants of a quadratic irrational,
    then  $\lim_{n \to \infty} f_n^2/F_n$ exists and is finite.
    \end{enumerate}
    \end{theorem}

    Statement 2 of Theorem \ref{ASY} motivates the following definitions.
    \newline
    \newline
    \noindent
    {\bf Definition:\/}
  Given a coloring $\cal C$ we define
  \begin{equation}
    \eta({\cal C\/})=\frac{f^2}{F},
  \end{equation}
  where $f$ is the fold count for $\cal C$ and $F$ is the
  number of triangles in the triangulation which $\cal C$ colors.
  We call $\eta({\cal C\/})$ the {\it Eisenstein Isoperimetric
    Ratio\/}
  of $\cal C$.
  \newline
  \newline
  {\bf Definition:\/}
  Given a quadratic irrational $\eta \in (0,1)$ let
  $$\eta(\zeta)=\lim_{n \to \infty} \eta(p_n +q_n \alpha),$$
  where $\{p_n/q_n\}$ is the sequence of
  continued fraction approximants of $\zeta$.  Theorem
  \ref{ASY} guarantees that this limit exists. We call
  $\eta(\zeta)$ the {\it Eisenstein Isoperimetric Ratio\/} of
  $\zeta$.
  \newline

In \S \ref{asy} we will show, among other calculations, that
  \begin{equation}
    \label{FIBO}
    \eta(\phi^{-1})=\phi^{6}, \hskip 30 pt
    \phi=\frac{\sqrt 5+1}{2}.
  \end{equation}
  Here $\phi$ is the golden ratio.
  Our proof of Statement 2 of Theorem \ref{ASY} will
  show more generally that $\eta(\zeta) \in \Q(\zeta)$,
  the quadratic field containing $\zeta$.
  See \S \ref{asy} for some examples.

  For comparison we prove the following easy result.
  \begin{theorem}
    \label{goodbound}
    Let $\cal C$ be any good coloring of any
    triangulation with degree sequence
    $6,...,6,2,2,2$.  Then
    $\eta({\cal C\/}) \geq 3$.
  \end{theorem}
  Theorem \ref{goodbound} indicates that some of the
  continued fraction colorings roughly minimize
  the fold count.  Here is a strong conjecture along
  these lines.

  \begin{conjecture}
\label{AI}
  Suppose $\{{\cal G\/}_n\}$ is any infinite sequence
  of distinct colorings of primitive
  sphere triangulations with degree sequence
  $6,...,6,2,2,2$.      Then $\liminf_{n \to \infty} \eta({\cal
    G\/}_n) \geq \phi^6$.
\end{conjecture}
Here {\it primitive\/} means that the triangulation has the form
${\cal T\/}(\beta)$ for a primitive Eisenstein integer $\beta$.
Equation \ref{FIBO} says that this conjecture
is sharp.  Since
$\phi^6=17.944...<18$, Theorem \ref{goodbound} says that the conjecture
is true up to a factor of $6$.

Sometimes there are good colorings which have lower fold count then
the corresponding continued fraction coloring. See Figure 9
in \S 3 for an example.  Here is one last conjecture about the general
situation.
\begin{conjecture}
  There is a constant $\Omega$ with the following
  property.  For any triangulation $\cal T$ with degree
  sequence $6,....,6,2,2,2$ there is a good coloring $\cal C$ of
  $\cal T$ with $\eta({\cal C\/})<\Omega$.
\end{conjecture}

    \subsection{Organization}

In \S 2, after a discussion of the slow Gauss map and its connection to
continued fractions, I will construct the
continued fraction colorings.
The building blocks  are what I call
{\it capped flowers\/}, and these in turn are made in layers
from cyclically arranged patterns of trapezoids which I call
{\it trapezoid necklaces\/}.
(Look again at Figure 1.)   I will explain how the set of
trapezoid necklaces is naturally the vertex set of the
infinite rooted binary tree (modified to have an extra vertex at the
bottom).  Taking a path in this tree defines the capped flower.
In \S 3 I will prove Theorem \ref{ASY},
establish Equation \ref{FIBO}, and give some evidence
for Conjecture \ref{AI}.
In \S 4 I will prove Theorem \ref{goodbound}.

\subsection{Acknowledgements}

I'd like to thank  Ethan Bove, Peter Doyle, Jeremy Kahn, Rick Kenyon, Curtis McMullen, Kasra Rafi, and
Peter Smillie for various conversations (sometimes going back some years)
on topics related to the material here.

 \newpage

\section{The Main Construction}

\subsection{The Slow Gauss Map}
\label{SG}

Given $r \in (0,1)$ we define
$\gamma(r)$ to be whichever of the two numbers
\begin{equation}
  \frac{1}{r}-1  \hskip 30 pt {\rm or\/} \hskip 30 pt
  \frac{1}{\frac{1}{r}-1}
\end{equation}
lies in $(0,1]$.
When $r>1/2$ the first choice works and
when $r<1/2$ the second choice works.
When $r=1/2$ both choices give $\gamma(1/2)=1$.
The map $\gamma: (0,1) \to (0,1]$ is called the
{\it slow Gauss map\/}.

Let $\cal S$ denote the set of rationals in $(0,1]$.
On ${\cal S\/}-\{1\}$ we have:
\begin{equation}
  \gamma\bigg(\frac{a}{a+b}\bigg)=\gamma\bigg(\frac{b}{a+b}\bigg)=\frac{a}{b}.
\end{equation}
We can think of $\cal S$ as the
infinite rooted binary tree (modified to have an extra bottom vertex).
We make this tree by joining each member of ${\cal S\/}-\{1\}$ to
its image under $\gamma$.   Figure 2 shows the beginning of this tree.

\begin{center}
  \resizebox{!}{1.5in}{\includegraphics{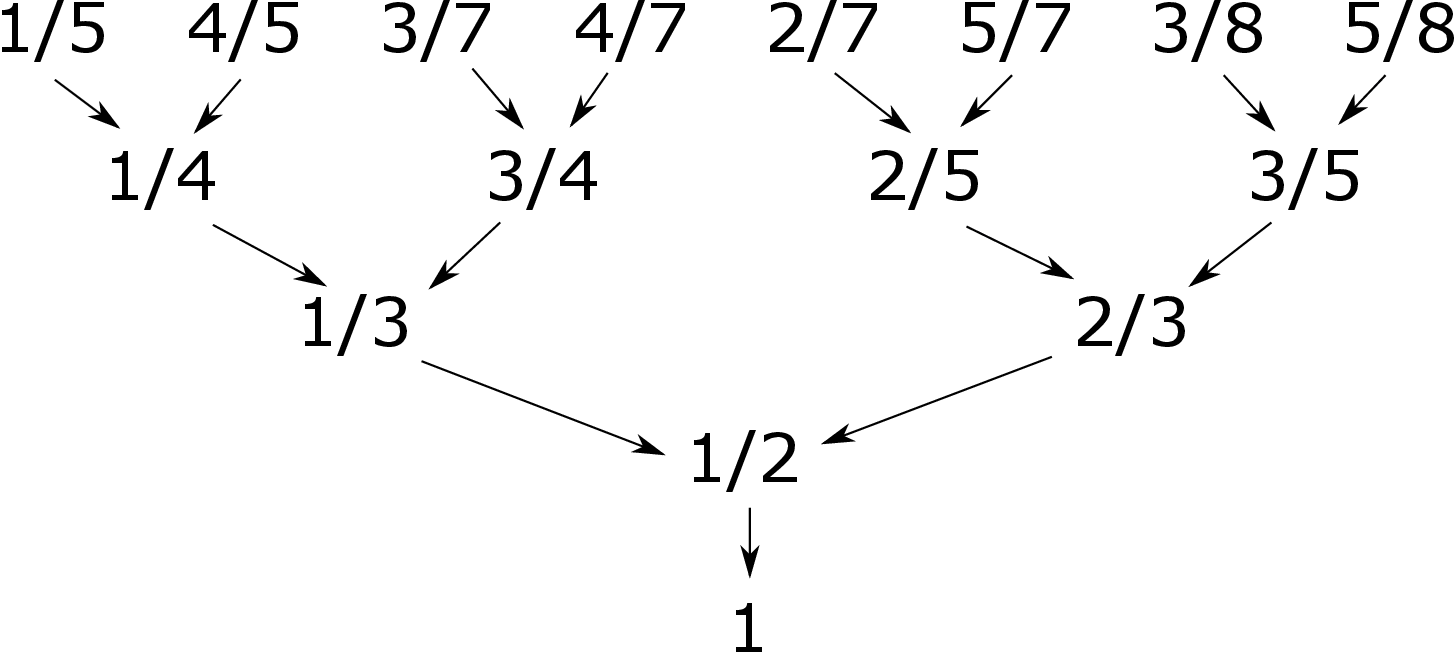}}
\newline
{\bf Figure 2:\/} The beginning of the tree of rationals.
\end{center}

The slow Gauss map is
connected to continued fractions.
The traditional Gauss map is
\begin{equation}
  \gamma^*(p/q)=(q/p)-{\rm floor\/}(q/p).
\end{equation}
For each $p/q$ there is some {\it comparison exponent\/} $k$ such
that $$\gamma^*(p/q)=\gamma^k(p/q).$$  In other
words, the (suitably) iterated slow Gauss map has the same action
as the traditional Gauss map; it just works more slowly.
The continued fraction expansion of $p/q$ is obtained by
recording the sequence of comparison exponents we see as
we iteratively apply $\gamma$ and $\gamma^*$ to $p/q$.

\subsection{Trapezoid Necklaces}

An {\it isosceles trapezoid\/} is a quadrilateral with two parallel sides
such that the other two sides are non-parallel but have the same length.
We call the longer parallel side the {\it top\/}, the
shorter parallel side the {\it bottom\/}, and the other two
sides the {\it diagonal sides\/}.
We allow the degenerate case of
an isosceles triangle.  In this case the bottom has length $0$.
An {\it Eisenstein trapezoid\/} is an isosceles trapezoid whose
edges lie in the
$1$-skeleton of $\cal T$, the planar equilateral triangulation whose
vertices are the Eisenstein integers. Figure 1 above and Figure 3 below
feature some Eisenstein trapezoids.

Up to symmetries of
$\cal T$, an Eisenstein trapezoid $X$ is characterized by
the pair $(a,b)$ where $a$ is the length of a diagonal
side of $X$ and $b$ is the length of the top.
We call $X$ {\it primitive\/} if
$a,b$ are relatively prime.
When $X$ is primitive, we define the
{\it aspect ratio\/} to be $a/b$.  The aspect
ratio determines the primitive Eisenstein trapezoid up to
symmetries of $\cal T$.   Thus, modulo symmetry,
the primitive Eisenstein trapezoids
are naturally in bijection with the set
$\cal S$ of rationals considered above.
The Eisenstein trapezoids in Figure 1 are
all primitive, and their aspect ratios are
variously $1/1$ and $1/2$ and $2/3$.
The Eisenstein trapezoids in Figure 3
have aspect ratio $3/5$.

\begin{center}
  \resizebox{!}{3.5in}{\includegraphics{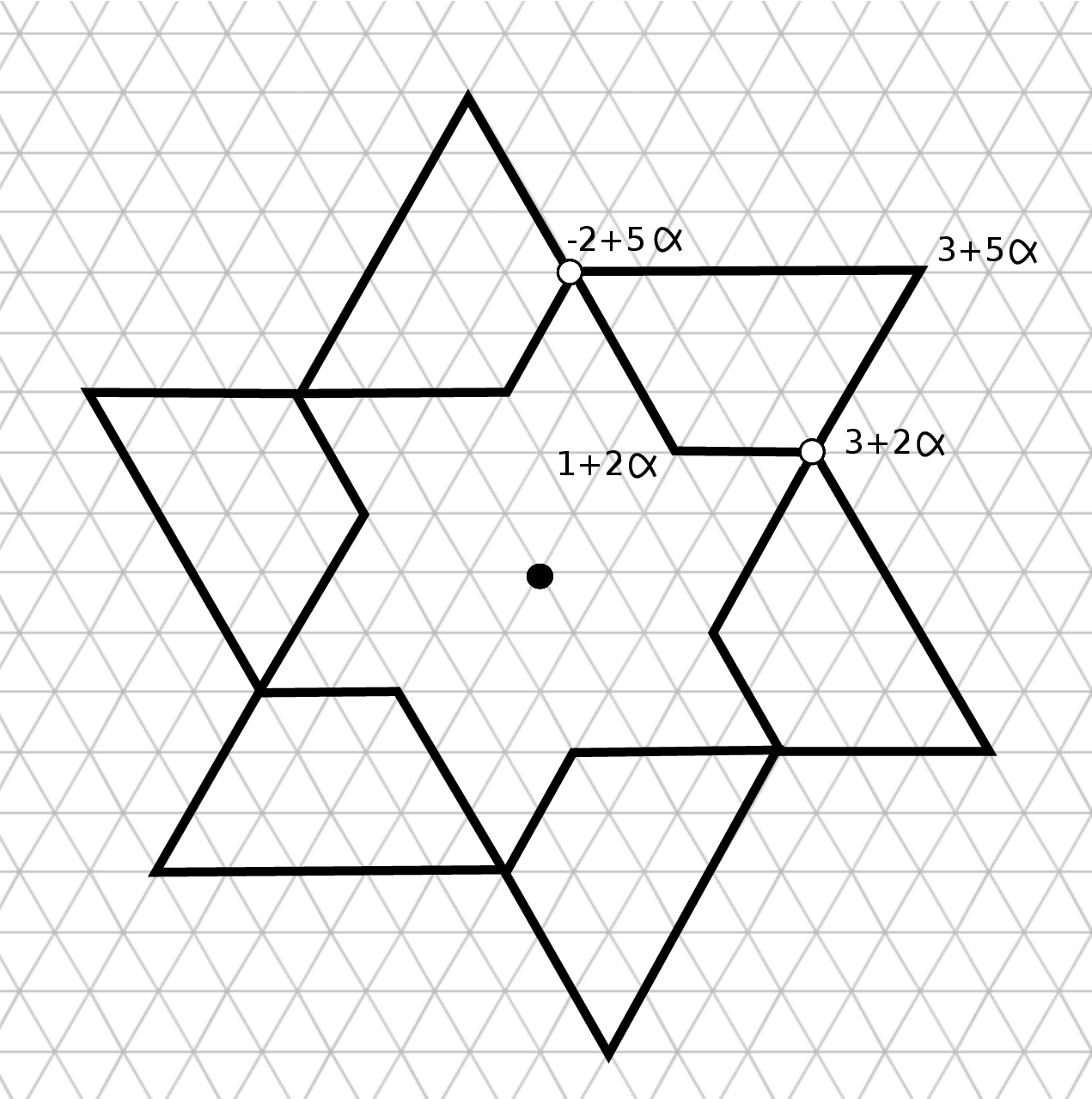}}
\newline
{\bf Figure 3:\/} A trapezoid necklace of aspect ratio $3/5$.
\end{center}

Figure 3 illustrates what we mean by a
{\it trapezoid necklace\/}. This is a union of
$6$ primitive Eisenstein trapezoids $X_1,...,X_6$ which has
the following properties:
\begin{itemize}
\item The trapezoids have pairwise disjoint interiors.
\item $X_i \cap X_{i+1}$ is a single point, a common vertex, for all $i$.
\item An order $6$ rotation $\rho$ of $\cal T$ has the action
  $\rho(X_i)=X_{i+1}$ for all $i$.
\end{itemize}
In this description the indices are taken mod $6$.
We define the {\it center\/} of the necklace to be
the fixed point of $\rho$.  When the center is $0$,
the map $\rho$ (or perhaps its inverse)
is multiplication by $\alpha$.
We define the aspect ratio of the necklace to be
the common aspect ratio of the $6$ individual
trapezoids.

Up to symmetry of $\cal T$, there exists a unique
Eisenstein necklace of aspect ratio $a/b \in \cal S$.
If we normalize the picture so that $0$ is the center,
then one of the trapezoids $X_1$ has vertices
$$a+b\alpha, \hskip 30 pt (a-b)+b\alpha,
\hskip 30 pt (2a-b)+ (b-a) \alpha, \hskip 30 pt a+(b-a)\alpha.$$
The intersection of $X_1$ and $X_2=\rho(X_1)$ is
the point $(a-b)+b\alpha$ because
$$\alpha \times \big(a+(b-a)\alpha)=(a-b)+b\alpha.$$
This little calculation uses the fact that
$\alpha^2=\alpha-1$.

\subsection{Empty Trapezoid Flowers}
\label{empty}

  Each trapezoid necklace $X$ defines a smaller
  trapezoid necklace $Y=\gamma(X)$ having the same center.
  The defining property is that the top side of each trapezoid $Y_i$ in
  $Y$ is a side of a trapezoid $X_j$ of $X$, and one of
  the diagonal sides of $Y_i$ is a side of one of the trapezoids
  of $X$ adjacent to $X_j$.  This awkward definition is very
  much like a written description of how to drink a glass of water.
  A demonstration says a thousand words.
  Figure 4 shows the trapezoid necklaces
  $$X \to \gamma(X) \to \gamma^2(X) \to \gamma^3(X)$$ alternately
colored black and white. Here $X$ is as in Figure 3.
The respective aspect ratios are given by
$3/5 \to 2/3 \to 1/2 \to 1/1$.
  We call the union of these necklaces the
  {\it empty\/} $3/5$-{\it flower\/}.

\begin{center}
  \resizebox{!}{5.75in}{\includegraphics{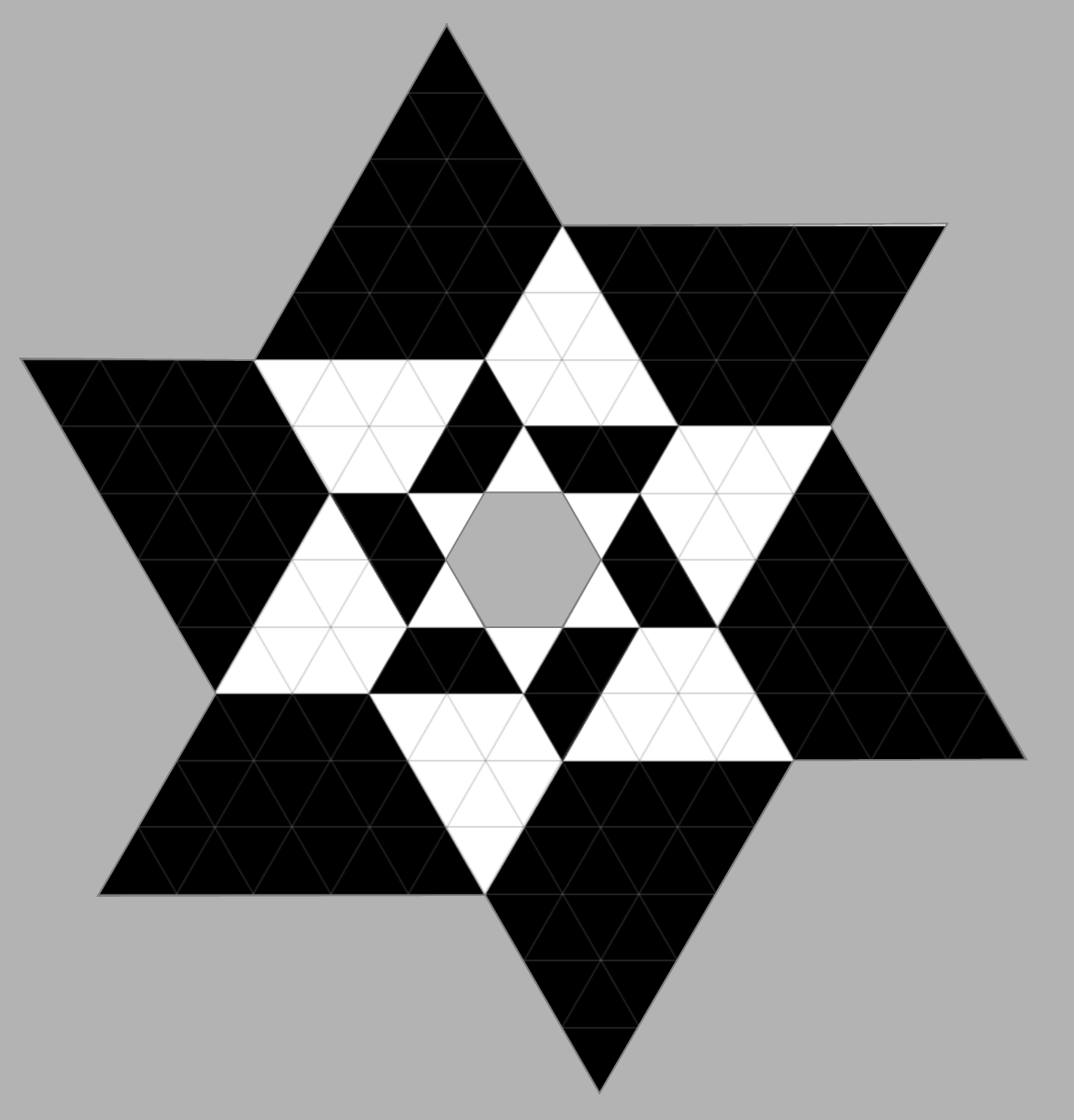}}
\newline
{\bf Figure 4:\/} The empty $3/5$-flower.
  \end{center}

  As is suggested by our notation, the action of
  $\gamma$ here mirrors the action of the
  slow Gauss map $\gamma$ defined above.
  That is, if $r$ is the
  aspect ratio of $X$ then $\gamma(r)$ is
  the aspect ratio of $\gamma(X)$.

  We discussed above how the slow Gauss map is
  related to continued fractions.  Here we
  continue the discussion.
  As we now illustrate, our construction also
  precisely implements the continued fraction expansion.

\begin{center}
  \resizebox{!}{4.8in}{\includegraphics{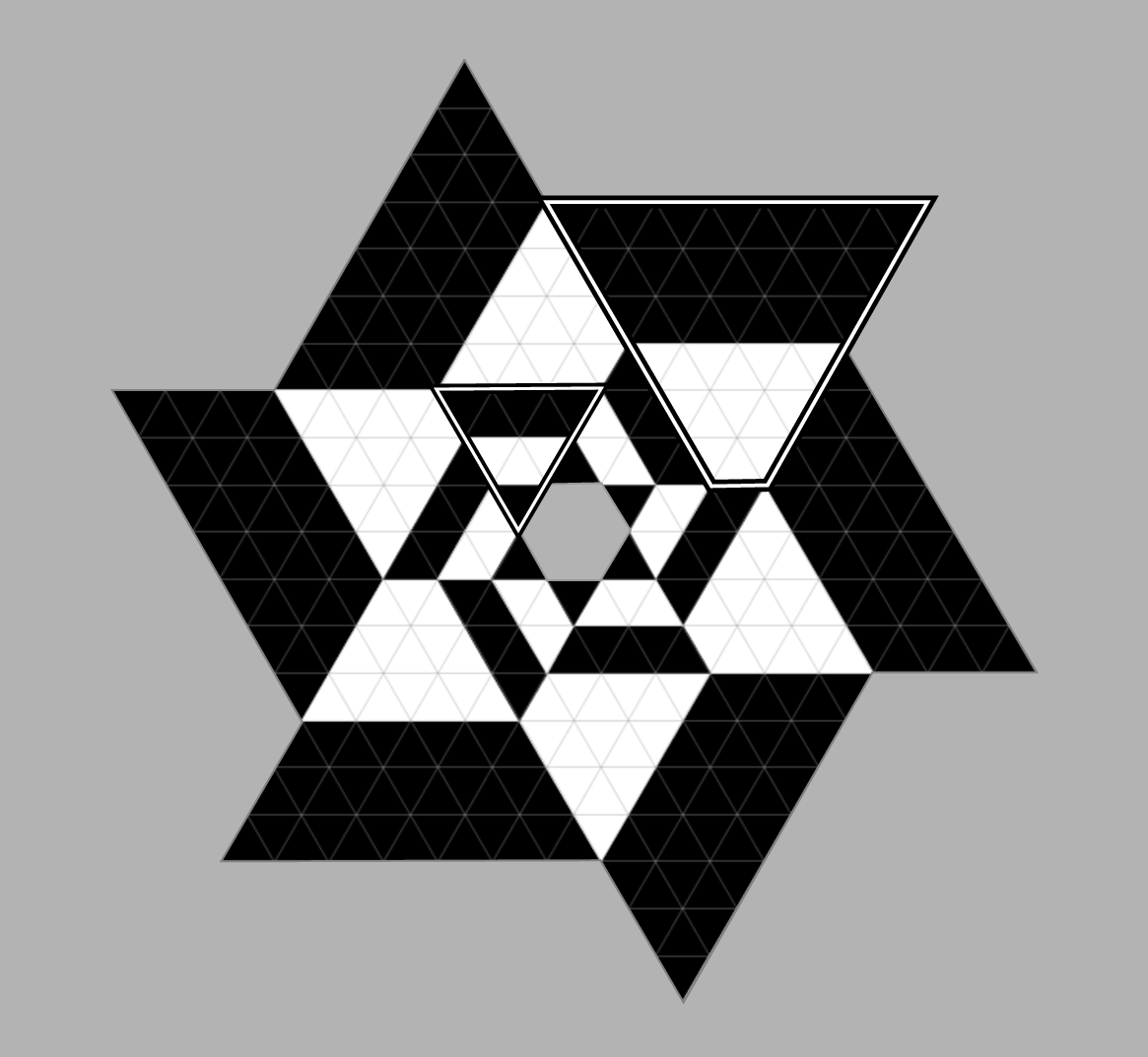}}
\newline
{\bf Figure 5:\/} The empty $3/7$-flower.
\end{center}

Figure 5 shows the empty $3/7$-flower, corresponding to the
tree path
$3/7 \to 3/4 \to 1/3 \to 1/2 \to 1/1$.  The specially-outlined
trapezoids (and their rotated images)
are the maximal trapezoids in the flower. They have
respectively $2$ and $3$ ``stripes''.
For comparison,
  $3/7$ has continued fraction $0:2:3$.  That is
$$\frac{3}{7}=0+\frac{1}{2+\frac{1}{3}}.$$

In general, the empty $p/q$-flower starts with a
$p/q$-necklace and then fills in the full
$\gamma$ orbit, alternately coloring the
necklaces black and white.  The innermost
necklace always has aspect ratio $1/1$.
Up to symmetries of $\cal T$, the empty
$p/q$-flower is unique.  One can read off
the continued fraction expansion of $p/q$ by
counting the stripes of the maximal
trapezoids.

\subsection{Filling and Capping the Flowers}

We fill an empty flower by coloring the remaining
$6$ triangles black and white in an alternating pattern.
Figure 6 shows this for the
$3/5$-flower. 
Up to symmetries of $\cal T$ (and
swapping the colors) there
 is a unique $p/q$-flower.  The empty flowers have
$6$-fold rotational symmetry but the (filled) flowers have
$3$-fold rotational symmetry.  We break a symmetry
to define the filling.

\begin{center}
  \resizebox{!}{4.4in}{\includegraphics{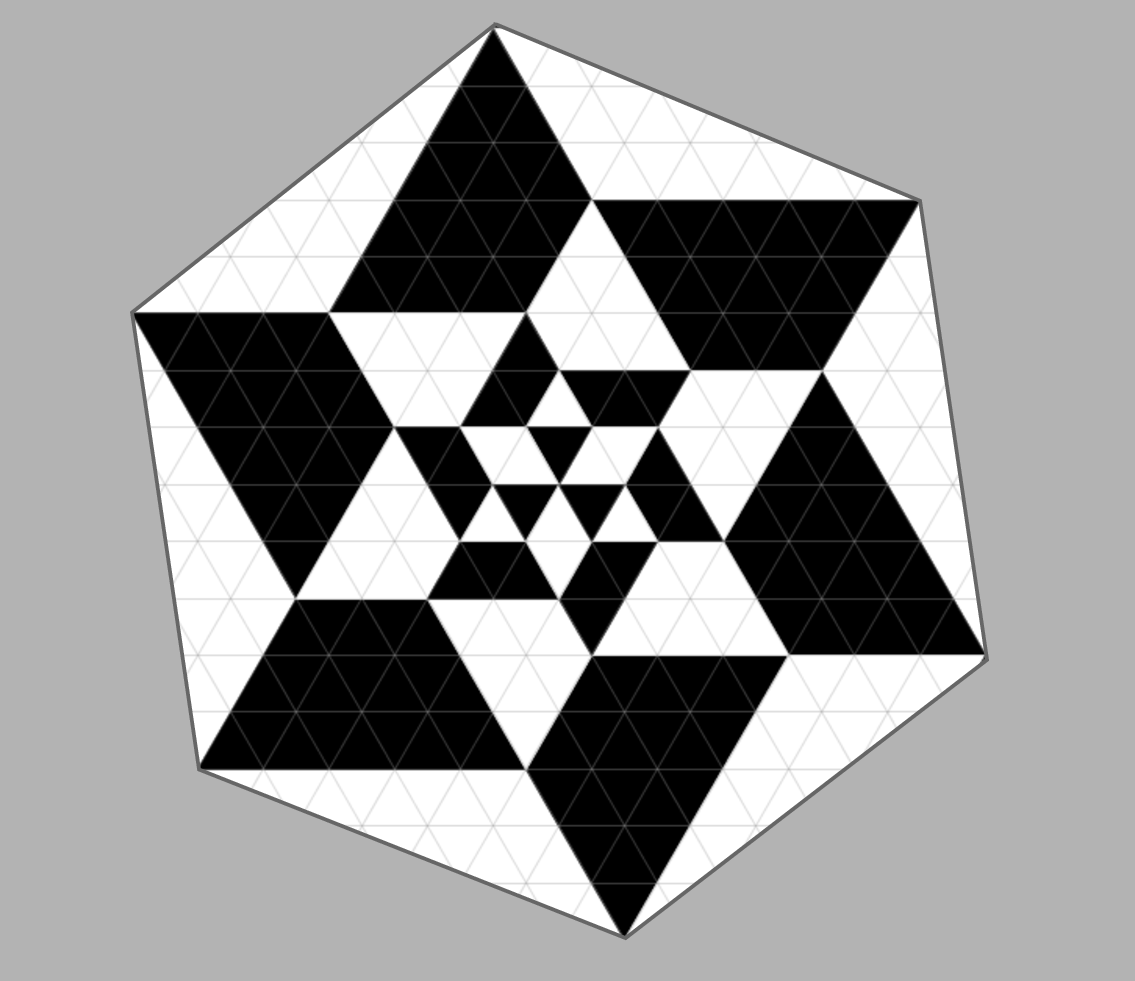}}
\newline
{\bf Figure 6:\/} The filled and capped $3/5$-flower
\end{center}

Figure 6 also shows what we mean by {\it capping\/} a flower.
We take the convex hull of the flower and color the complementary
triangles the color opposite the color of the outer necklace.
The capped flowers, which are regular hexagons with
Eisenstein integer vertices, are the building blocks of our
colorings.  When two translation-equivalent
capped flowers meet
along a common boundary edge, the triangular regions
merge to become an {\it Eisenstein parallelogram\/} -- i.e.,
one whose boundary lies in
the $1$-skeleton of $\cal T$.

\subsection{Defining the Colorings}

Consider the hexagonal tiling of the plane by translates of the
capped $a/b$-flower. Because the capped flowers have
$3$-fold rotational symmetry, the resulting planar coloring
is invariant under the group $G_{\beta}$ generated by order $3$
reflections in the vertices and centers of the hexagons.
The quotient of this planar coloring by
$G_{\beta}$ is ${\cal C\/}(\beta)$.
By construction, ${\cal C\/}(\beta)$ is good.

\begin{center}
  \resizebox{!}{4in}{\includegraphics{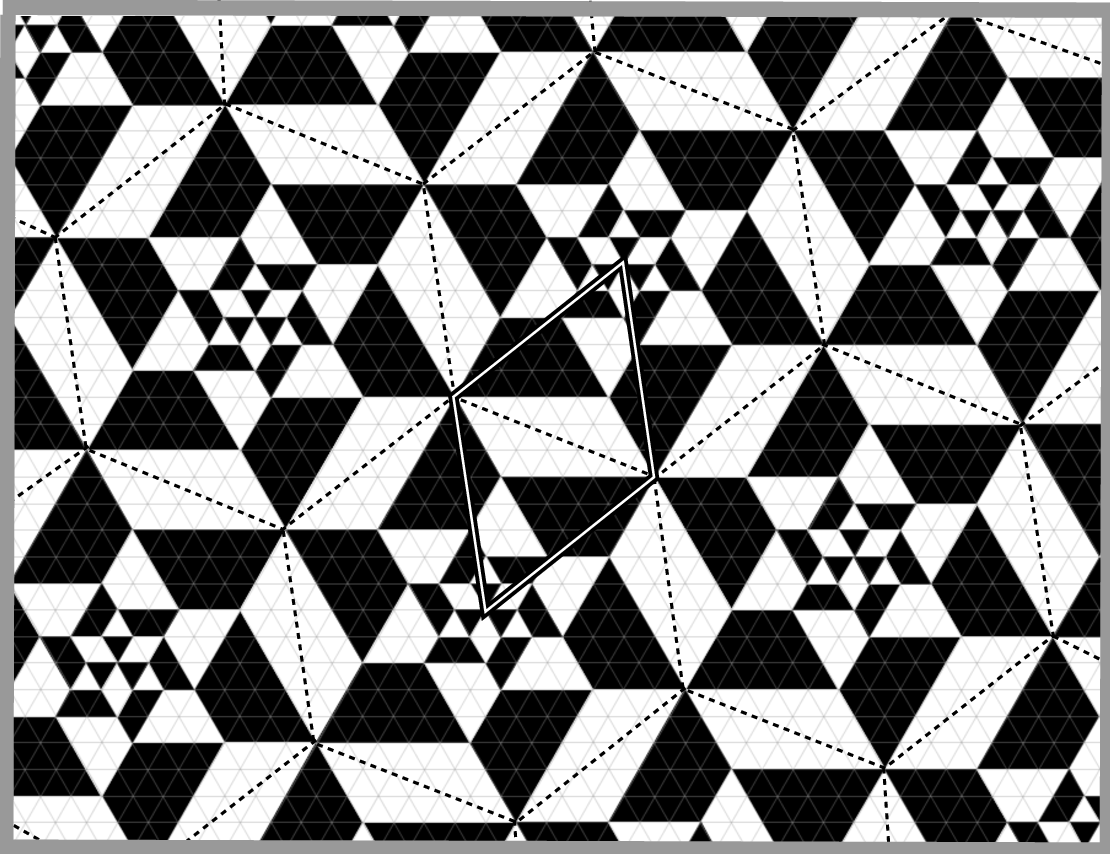}}
\newline
{\bf Figure 7:\/} The universal cover of ${\cal C\/}(3+5\alpha)$.
\end{center}

Figure 7 shows the construction for $\beta=3 + 5 \alpha$.
The region bounded by the big central rhombus is
a fundamental domain for the action of $G_{\beta}$.
The colorings exhibit a lot of variety.
Figure 8 below shows ${\cal C\/}(a+13 \alpha)$ for
$a=1,3,5,7,9,11$.   It is worth noting that
for parameters like $1/13$ the fold count is
linear in the number of triangles.

\begin{center}
  \resizebox{!}{5.6in}{\includegraphics{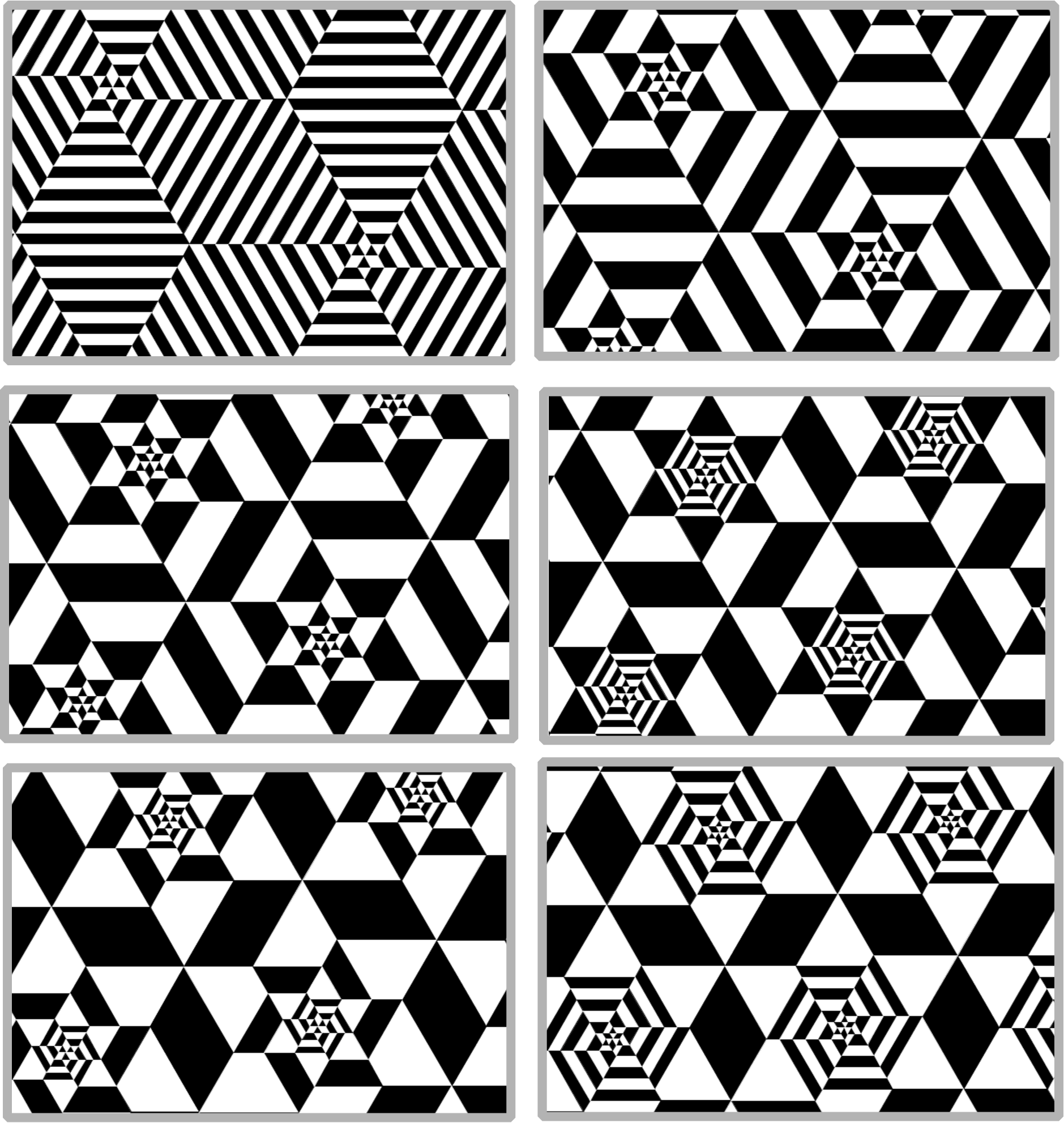}}
\newline
{\bf Figure 8:\/} The orbifold covers of ${\cal C\/}(a+13\alpha)$
for $a=1,3,5,7,9,11.$.
\end{center}

\newpage

\section{Properties of the Colorings}

\subsection{Zero Degree}
\label{even}

Let $T$ denote the regular tetrahedron.
Call a triangulation of the sphere {\it even\/} if it
has all even degrees. In this section we prove
that any good coloring of an even triangulation $\Sigma$
has the same total number of black and white triangles.
This is equivalent to the statement that the
associated map $f: \Sigma \to T$ has topological degree $0$.
I learned this proof from Kasra Rafi.

We equip $\R^2$ with the
unit equilateral triangulation.
Let $\pi: \R^2 \to T$ be the branched covering
which maps triangles to triangles in an
isometric way.  A good triangulation has the property
that, around each vertex, the number of
black triangles is congruent mod $3$ to the
number of white triangles.  In the case of
an even triangulation, a good coloring
has the stronger property that around each vertex the number of
black triangles is congruent mod $6$ to the
number of white triangles.
But this means that the map $f: \Sigma \to T$
lifts to a map $\widetilde f: \Sigma \to \R^2$ such that
$f=\pi \circ \widetilde f$.  The map $\widetilde f$ is
homotopic to a point and hence so is $f$.   But then
$f$ has topological degree $0$.

\subsection{Asymptotic Properties}

In this section we will prove Theorem \ref{ASY}.
Recall that $\cal S$ is the set of rational numbers in $(0,1]$.
Now suppose that $\{p_n/q_n\}$ is an infinite sequence of
elements of $\cal S$ having an irrational limit $\zeta$.
\newline
\newline
{\bf Statement 1:\/}
If we have an
isosceles trapezoid based on the pair
$(a,b)$ then the length-to-area ratio is
at most $100/a$.  This is conservative estimate.
We fix
$D$ for the whole discussion and show that
$ f_n/F_n <200/D$ once
$n$ is large enough.
Letting $D \to \infty$ we see that
$f_n/F_n \to 0$.

Let $G(p_n/q_n)$ denote the sequence of rationals
we get by starting with $p_n/q_n$, iteratively applying the slow
Gauss map, then listing the terms in reverse order.
For instance,
$$G(3/7)=\frac{1}{1}, \hskip 10 pt
\frac{1}{2}, \hskip 10 pt \frac{1}{3}, \hskip 10 pt
\frac{3}{4}, \hskip 10 pt \frac{3}{7}.$$
We call this the $G$-{\it sequence\/}.

Let $p'_n/q'_n$
be the term furthest along in the $G$-sequence
with $p_n' \leq D$.   Let $R_n$ be the radius
of the flower based on $p_n/q_n$ and likewise define
$R'_n$.  

\begin{lemma}
  \label{fat}
  $R_n'/R_n \to 0$ as $n \to \infty$.
\end{lemma}

\startproof
There is a uniform constant $C>0$ such that
\begin{equation}
  \label{unif}
  \frac{1}{C} < \frac{R_n}{p_n+q_n}, \frac{R_n'}{p_n'+q_n'}<C.
\end{equation}
So, it suffices to prove that
$\rho_n:=(p_n'+q_n')/(p_n+q_n) \to 0$ as $n \to \infty$.

Say that a {\it big step\/} happens in the $G$-sequence
whenever the numerator for consecutive terms changes.
Whenever a big step happens,
the number $a/b$ precedes the number $b/(a+b)$.
The ratio of the sums is $(a+b)/(a+2b)<2/3$.
So, if $m_n$ big steps are taken between $p'_n/q'_n$
and $p_n/q_n$ then $\rho<(2/3)^{m_n}$.

Let $M_n$ be the total number of big steps taken
in $G(p_n/q_n)$.   Note that $m_n \geq M_n-D$.
By continuity, 
the Gauss map for $p_n/q_n$ matches the Gauss map
for the limit $\zeta$ for a number of initial steps that
tends to $\infty$ as $n \to \infty$.
Hence $M_n \to \infty$ as $n \to \infty$.
Hence $m_n \to \infty$ as $n \to \infty$.
Hence $\rho_n \to 0$.
\endproof

\noindent
{\bf Remark:\/} I am grateful to Emma R\"uter for pointing
out a flaw in the version of Lemma \ref{fat} that I had in
an earlier version of this paper.
\newline

The little flower based on $p_n'/q_n'$ is a subset of the big
flower based on $p_n/q_n$.  We have
\begin{equation}
  \label{outside}
  f_n'<100F_n', \hskip 30 pt \frac{f_n-f_n'}{F_n-F_n'} \leq \frac{100}{D}.
\end{equation}
The first bound holds in any flower. The second bound
comes from the fact that we are only considering trapezoids
outside the small flower, and these all have length-to-area ratio
at most $100/D$.

For any $\epsilon>0$ we can take $n$ so large that
$R_n'/R_n<\epsilon$.  But then, by scaling, we have
$F_n'/F_n<\epsilon^2$.    Choosing $n$ sufficiently large, so
as to make the final inequality true, we have
$$
  \frac{f_n}{F_n} =
  \frac{(f_n-f_n')+f_n'}{F_n} \leq  
  \frac{(f_n-f_n')+100F_n'}{F_n} \leq  $$
  \begin{equation}
    \frac{(f_n-f_n')}{F_n} +100\epsilon^2 <
    \frac{f_n-f'_n}{F_n-F'_n} + 100 \epsilon^2  \leq
  \frac{100}{D} + 100\epsilon^2<\frac{200}{D}.
\end{equation}
This proves Statement 1.
\newline
\newline
{\bf Statement 2:\/}
Recall that the quadratic irrational limit $\zeta=\lim p_n/q_n$ has an
eventually periodic continued fraction expansion and
(equivalently) is pre-periodic under the action of the slow Gauss map.

We think of ${\cal T\/}(p_n+q_n \alpha)$ as
a coloring of a doubled equilateral triangle.  We then scale the metric
by a factor of $q_n^{-1}$.  The corresponding sequence of doubled equilateral
triangles converges (say, in the Gromov-Hausdorff topology)
to another doubled equilateral triangle, namely
$\C/G_{\beta}$ where
$\beta=\zeta + \alpha.$
Again, $G_{\beta}$ is the group generated by order $3$
rotations in points of $\beta {\cal E\/}$.

The rescaled colorings also converge.  The limiting coloring
is made from ``infinite flowers'' and parallelograms.  The infinite flowers
are limits of the rescaled flowers associated to the Eisenstein
integers $a_n+b_n\alpha$.   Scaling appropriately, we arrange that
one of the outer trapezoids in the infinite flower $F$
has aspect ratio $\zeta$,  top right vertex $\beta$
and horizontal top and bottom.   The remaining trapezoids have
aspect ratios which we find by applying the slow Gauss map iteratively
to $\zeta$.   These necklaces 
converge to the origin.

Let $F[k]$ denote the infinite flower obtained by trimming off the
outermost $k$ trapezoid necklaces of $F$.  Given the preperiodicity
of the slow Gauss map on $\zeta$,
there are integers $m,n$ and
some $\lambda \in (0,1)$ such that
\begin{equation}
  F[m+kn]= \lambda^k F[m], \hskip 30 pt k=1,2,3,...
\end{equation}
In other words, if we strip off the outer $m$ layers of $F$ then the
resulting infinite flower is self-similar with ``period'' $n$ and similarity factor $\lambda$.

The area of a unit equilateral triangle is $\sqrt 3/4$ and the
side length is $1$.  Let us redefine the Eisenstein isoperimetric
ratio to be the quantity
\begin{equation}
  \label{Iso2}
  \frac{{\rm perimeter\/}^2}{{\rm area\/}}\times \frac{4}{\sqrt 3}.
\end{equation}
This new definition coincides with the old definition for colorings
based on unit equilateral triangulations.  Also, the new definition
is scale invariant.  Thus, to prove Statement 2 of
Theorem \ref{ASY} we just need to show that
this quantity is finite and well-defined for our limiting coloring.
This suffices because the invariant then varies continuously.

We first compute the area.  This is given by
$${\rm area\/}=C_1 + C_2 \sum_{k=1}^{\infty} \lambda^{2k} = C_1 + 
\frac{C_2}{1-\lambda^2}.$$
Here $C_1$ is a constant that depends on the $m$ outer layers of
the infinite flower and on the triangular regions defining the cap
of the flower $F$.  The constant $C_2$ is determined by the
outer $n$ layers of $F[m]$.

Now we compute the perimeter.  This is given by
$$
{\rm perimeter\/} =C_3 + C_4 \sum_{k=1}^{\infty} \lambda^k=C_3 +
\frac{C_4}{1-\lambda}$$
The constants $C_3,C_4$ have a similar dependence as
$C_1,C_2$.
So, the quantity in Equation \ref{Iso2} exists and is finite.
This completes the proof of Statement 2.
\newline
\newline
{\bf Remarks:\/} \newline
(1)
The quantities $C_1 \sqrt 3$ and $C_2 \sqrt 3$ both belong
to $\Q(\zeta)$.  Likewise $C_3,C_4, \lambda$ also all belong to
$\Q(\zeta)$.
Therefore, the limiting value $\eta(\zeta)$ lies in
$\Q(\zeta)$.  I will give some calculations below.
\newline
(2) Our insistence that $\{p_n/q_n\}$ be the sequence
of continued fraction approximants of $\zeta$ is more restrictive than
need be.  The limiting argument works as long as
$\{p_n/q_n\}$ limits to $\zeta$ and the corresponding
continued fraction expansions are uniformly bounded.
Thus, the Eisenstein isoperimetric ratio is well-defined
for any irrational in $(0,1)$ with bounded continued
fraction approximation.

\subsection{Some Calculations}
\label{asy}

Let us first establish Equation \ref{FIBO}.
Rather than take the general approach from the
last section we work more concretely with
the sequence $\{a_n/a_{n+1}\}$ of continued
fraction approximants to $\phi^{-1}$. Here
$$a_1,a_2,a_3,a_4,a_5,...1,1,2,3,5,...$$ is the
sequence of Fibonacci numbers.  We will
use the well-known asymptotic formula
\begin{equation}
    \label{approx}
    a_n \sim \frac{\phi^n}{\sqrt 5}.
  \end{equation}

  For $n \geq 2$ let
  $f_n$ denote the fold count of ${\cal C\/}(a_n/a_{n+1})$ and
  let $F_n$ denote the number of triangles.
  From our program we compute the following table:
  \begin{center}
    \begin{tabular}  { c c c  c c c}
     $n$ & $a_n$ & $a_{n+1}$ & $f_n$ & $F_n$ & $\phi^{-6} f_n^2/F_n$
      \\
      \\
  2& 1& 2 & 13 & 14 & .672718... \\
   3& 2&3 & 23 &38 & .775795...\\
   4& 3&5 & 39 & 98 & .864923...\\
  5&  5&8 & 65 & 258 & .912601...\\
  6&  8&13 & 107 & 674 & .946633...\\
  \end{tabular}
  \end{center}
  After a little trial and error we find the following formulas
\begin{equation}
  f_n=-1+2 \sum_{k=1}^{n+2} a_k, \hskip 30 pt
  F_n=2 + 4 \sum_{k=1}^{n} a_k a_{k+1}.
\end{equation}
These formulas give the same answers as in the table above for $n=2,3,4,5,6$.
An easy inductive proof, which we omit, establishes them.
Using the approxiation in Equation \ref{approx}, and the familiar
formula for the partial sums of a geometric series, and the
fact that $1+\phi=\phi^2$, we find that
$$f_n^2 \approx \frac{4}{5} \times \phi^{2n+8}, \hskip 30 pt
F_n \approx \frac{4}{5} \times \phi^{2n+2}.$$
Here $\approx$ means equal up to a uniformly bounded error.
Dividing the one equation by the other gives
Equation \ref{FIBO}.

Now we describe some additional (nonrigorous) calculations we
did in Mathematica [{\bf Wo\/}].   Starting with a quadratic
irrational $\zeta$ we do the following.
\begin{enumerate}
  \item We take a close continued fraction approximation $p/q$
to $\zeta$.  For our calculations we used the command
{\bf Rationalize[$\zeta$,Power[10,-1000]]\/}.  This gives us a
continued
fraction approximant which agrees with $\zeta$ up to $1000$ decimal
places.  (We needed roughly this much accuracy in a few cases.)
\item We compute
  the continued fraction expansion of $\eta(p + q \alpha)$.
  \item We extract
what appears to be the start of a preperiodic continued fraction
expansion for a quadratic irrational and we then reconstruct this
quadratic irrational from the (implied) preperiodic continued
fraction expansion.  (In all cases, the fractional part of the
expression was strictly periodic.)
\item After finding what appears to be the period of the continued
  fraction expansion, we reconstruct the quadratic irrational
  which has this continued fraction expansion.  We guess that this is
  probably
  $\eta(\zeta)$.
\end{enumerate}
Here we show some examples.  Let
$\zeta_n=\sqrt{n}-{\rm floor\/}(\sqrt n).$  We have
\begin{equation}
  \begin{matrix}
    \eta(\zeta_2) &= & \frac{75 + 53 \sqrt 2}{7}  & &
    \eta(\zeta_3)&=& \frac{132+72 \sqrt 3}{13} \cr
    \cr
    \eta(\zeta_5)&=& \frac{321+137 \sqrt 5}{19} & &
    \eta(\zeta_6)&=& \frac{27+9 \sqrt 6}{2} \cr
    \cr
    \eta(\zeta_7)&=& \frac{3100+856 \sqrt 7}{259} & &
    \eta(\zeta_8)&=& \frac{1569+370 \sqrt 8}{98} 
  \end{matrix}
  \end{equation}
One could likely adapt the above proof of
Statement 2 of  Theorem \ref{ASY} to rigorously justify these calculations.

  \subsection{Evidence for the Isoperimetric Inequality}

  Most of our evidence for Conjecture \ref{AI} comes from
calculations within
the famiy of continued fraction colorings.
Given primitive Eisenstein integers
$$\beta=a+ b\alpha, \hskip 30 pt \beta' = a'+ b' \alpha,$$
we
write $\beta \preceq \beta'$ if $b \leq b'$.
It seems that
\begin{equation}
  \label{IE}
      \eta({\cal C\/}(\beta)) < \eta({\cal C\/}(\beta'))
      \end{equation}
    whenever $\beta \preceq \beta'$ and $\beta$ is a Fibonacci
    Eisenstein integer and $\beta' \not = \beta$.
    
    I checked Equation \ref{IE} when $\beta$ is any of
    the first $10$ examples and for all $|\beta'|<500$.  I think
    that this is very strong evidence.
    The proof of Equation \ref{IE} in general should be purely
    a matter of number theory.
    I haven't yet looked for a proof.
    Equation \ref{IE} would establish the Eisenstein Isoperimetric
    Inequality for the continued fraction colorings. These, of course,
    are only a tiny fraction of the good colorings.

    The only other evidence I have is the fact that
    ${\cal C\/}(1+2\alpha)$ and ${\cal C\/}(2+3\alpha)$ seem
    to minimize the fold counts for good colorings of
    ${\cal T\/}(1+2\alpha)$ and ${\cal T\/}(2 + 3 \alpha)$
    respectively.  My program is not good enough to
    convincingly check this for larger Fibonacci triangulations.
    Combining this scant evidence with the very good
    evidence for Equation \ref{IE}, I make the following conjecture.
    
    \begin{conjecture}
      Suppose that $\beta$ is a Fibonacci Eisenstein integer and
      that $\beta' \in (0,1)$ is a different Eisenstein integer such that
  $\beta \preceq \beta'$.  Let
  ${\cal C\/}={\cal C\/}(\beta)$. Let
  ${\cal C\/}'$ be an arbitrary
  good coloring of ${\cal T\/}(\beta')$.
  Then  we have $\eta({\cal C\/})<\eta({\cal C\/}')$.
    \end{conjecture}
  This conjecture would combine with Equation \ref{FIBO} to prove
  Conjecture \ref{AI}.

It is worth pointing out that that the continued fraction colorings
are far from the best in many cases.  Figure 9 compares the
continued fraction ${\cal C\/}(1+5\alpha)$ (left)  with another coloring
of ${\cal T\/}(1+5 \alpha)$ (right) which seems to be the one which
minimizes the fold count.

\begin{center}
  \resizebox{!}{3in}{\includegraphics{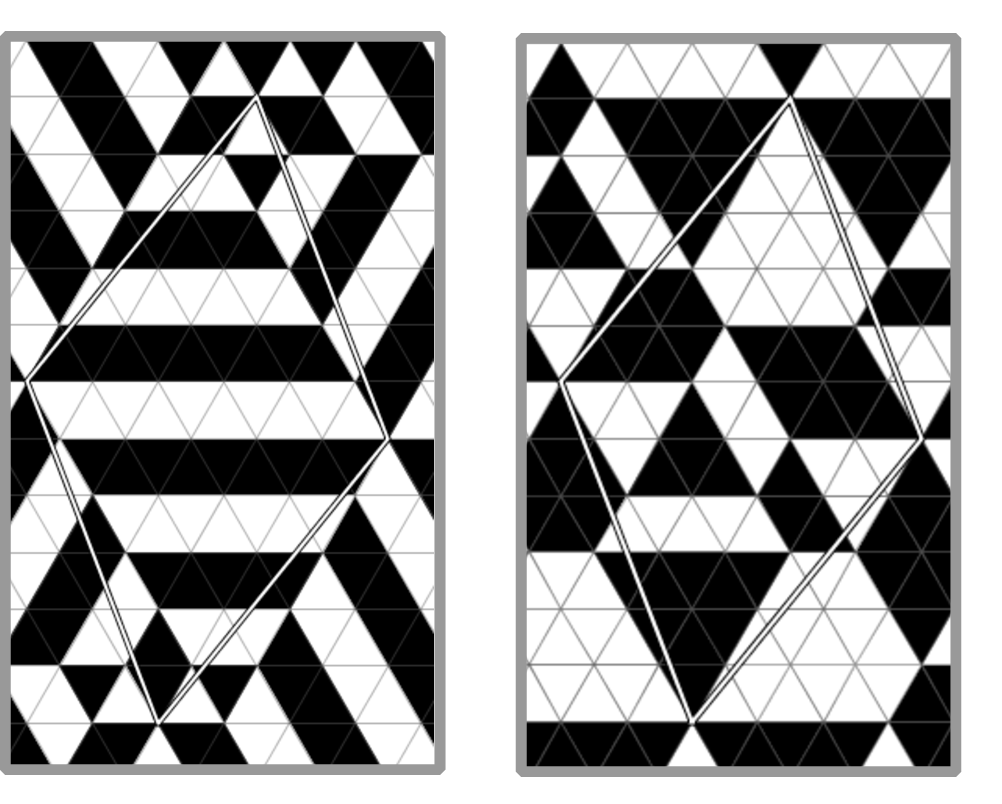}}
\newline
{\bf Figure 9:\/} Two good colorings of
${\cal T\/}(1+5\alpha)$.
  \end{center}

\newpage

\section{A Bound on the Isoperimetric Ratio}

In this chapter we prove Theorem \ref{goodbound}.
We first establish an isoperimetric inequality for a
special kind of polygon and then we give the
main argument.

\subsection{Special Hexagons}

Say that a {\it special hexagon\/} is a convex hexagon
whose interior angles are all $2\pi/3$.
The regular hexagon is an example, and all other
special hexagons are obtained from the regular
one by pushing the sides in and out parallel to
themselves.  Our first result is a version of the
isoperimetric inequality for special polygons.
This result is undoubtedly well-known.

\begin{lemma}
  \label{IS}
  The length and area of a special hexagon satisfy the inequality
  \begin{equation}
  \frac{{\rm length\/}^2}{{\rm area\/}} \geq 8 \sqrt 3.
  \end{equation}
\end{lemma}

\startproof
 It is an elementary exercise, and it is also discussed at length in
  [{\bf T\/}], that the area of such a hexagon is a
quadratic function of the side lengths
$\ell_1,...,\ell_6$.    By symmetry, the desired function
is symmetric in the arguments.  Therefore, we have
$${\rm area\/} = a(\ell_1+...+\ell_6)^2 + b (\ell_1^2 + ... + \ell_6^2).$$
Considering the unit regular hexagon and the unit equilateral
triangle, which is a limiting case, we get the relations
$$36 a + 6 b = \frac{3 \sqrt 3}{2}, \hskip 30 pt
9a + 3b = \frac{\sqrt 3}{4}.$$
Solving these equations leads to the formula
\begin{equation}
  {\rm area\/} = \frac{1}{6 \sqrt 3}(\ell_1+...+\ell_6)^2 -\frac{1}{4 \sqrt 3} (\ell_1^2 + ... + \ell_6^2).
\end{equation}
If we hold the perimiter fixed, we maximize the area by minimizing the
sum of the squares of the lengths.  This happens when all lengths are
equal, and the ratio of interest is scale-invariant.
Plugging in $\ell_1=...=\ell_6=1$ and calculating, we get
the advertised inequality.
\endproof

\subsection{The Main Argument}

We define an {\it Eisenstein polygon\/} to be a convex polygon whose
sides are contained in the $1$-skeleton of the unit equilateral triangulation.
Up to scaling, Eisenstein polygons are all limits of the
special hexagons considered in the last section. Let $f(P)$ denote the
perimeter of an Eisenstein polygon and let $F(p)$ denote the number of
triangles it contains.  As an immediate consequence of Lemma \ref{IS}, we have
\begin{equation}
  \label{IS2}
  \frac{f^2(P)}{F(P)} \geq 6.
\end{equation}

Let $\cal C$ be a good coloring of
some triangulation $\cal T$ with
degree sequence $6,...,6,2,2,2$.
Let $f$ denote the fold count
and let $F$ be the number of faces.
Our goal is to show that
$\eta({\cal C\/})=f^2/F \geq 3$.
Note that $\cal C$ partitions
$\cal T$ into a finite number of
monochrome regions with connected
interior.

\begin{lemma}
  \label{eis}
  Each monochrome region is isometric to
  an Eisenstein polygon.
\end{lemma}

\startproof
Let $P$ be such a region.  Without loss of generality assume that
$P$ is colored white.
To make the argument clearer, we shave off the
outer $\epsilon$ of $P$, so that its boundary 
does not self-intersect.
Let $P_{\epsilon}$ be this slightly smaller polygon.
Since both colors occur at the triangles touching
each degree $2$ vertex, we see that $P_{\epsilon}$
does not contain any degree $2$ vertices.
Hence $P_{\epsilon}$ is locally Eucidean.

Consider $\partial P_{\epsilon}$.
Each vertex $v_{\epsilon}$ of $\partial P_{\epsilon}$ is within
$\epsilon$ of a unique vertex $v$ of $P$.
Because $P$ comes from a good coloring, $v$
  has at most $3$ white
  triangles around it.   Hence $P_{\epsilon}$ is locally
  convex near each boundary component.

  Consider a lift $\widetilde P_{\epsilon}$ of
  $P_{\epsilon}$ to $\C$, the orbifold universal
  cover of our flat cone sphere.
  This lift is also locally convex near its
  boundary components.  If $\widetilde P_{\epsilon}$
  is unbounded then the local convexity near the boundary
  forces it to be an infinite strip.  But then $P_{\epsilon}$
  is an annulus with geodesic boundary. This is impossible
  in our flat cone sphere.    Hence $\widetilde P_{\epsilon}$ is
  bounded.  The local convexity forces $\widetilde P_{\epsilon}$
  to have just
  one boundary component.  If the map $\widetilde P_{\epsilon} \to
  P_{\epsilon}$ is not a global isometry then $P_{\epsilon}$
  contains some degree $2$ vertices, a contradiction. We conclude
      that $\widetilde P_{\epsilon}$ is a convex polygon and
    the map $\widetilde P_{\epsilon} \to P_{\epsilon}$ is an
    isometry.  Hence $P_{\epsilon}$ is a convex polygon.
    
  Letting $\epsilon \to 0$ we see that $P$ is also a convex polygon.
  But then $P$ is also an Eisenstein polygon.
  \endproof

Let $P_1,...,P_k$ be the white connected 
regions.  By Lemma \ref{eis},  each $P_j$ is an
Eisenstein polygon.
Let $f_j=f(P_j)$ and $F_j=F(P_j)$.
We have
\begin{equation}
  \min_j \frac{f_j^2}{F_j} \geq 6, \hskip 30 pt
  f=f_1+...+f_k, \hskip 30 pt
  \frac{F}{2}=F_1+...+F_k.
\end{equation}

The middle equality comes from the fact that
each black-white interface in particular lies
in the boundary of some $P_j$, and that
$\partial P_i$ and $\partial P_j$ are disjoint except perhaps for
vertex intersections when $i \not =j$.
Hence
\begin{equation}
 \frac{2f^2}{F} =\frac{(f_1+...+f_k)^2}{F/2}\geq \frac{f_1^2+...+f_k^2}{F/2}= \frac{f_1^2}{F_1}
 \oplus ... \oplus \frac{f_k^2}{F_k} \geq 6.
\end{equation}
Here $\oplus$ denotes the Farey sum:  We add the numerators
and we add the denominators.  As is well known, the Farey sum of
some fractions lies between the minimum and the maximum of
the fractions.

\section{References}

\noindent
[{\bf T\/}],  W. P. Thurston, {\it Shapes of polyhedra and
  triangulations of the sphere\/},  Geometry and Topology monographs
1: pp 511--549 (1998)
\newline
\newline
[{\bf W\/}], D. B. West, {\it Introduction to Graph Theory, 2nd Ed.\/},
  Prentice-Hall (2000)
  \newline
  \newline
[{\bf Wo\/}] S. Wolfram, {\it Mathematica\/} (2020) wolfram.com/mathematica.

\end{document}